\documentclass[aap,MSNbibl,citesort,dvips]{arximspdf}

\doi{10.1214/11-AAP772}
\volume{22}
\issue{1}
\pubyear{2012}
\firstpage{337}
\lastpage{362}

\makeatletter
\newcommand{\eqref}[1]{(\ref{#1})}
\newcommand{\E}{\mathbb{E}}
\newcommand{\R}{\mathbb{R}}
\newcommand{\N}{\mathbb{N}}

\newcommand{\Z}{\mathbb{Z}}
\newcommand{\VVV}{\mathbb{T}}

\renewcommand{\P}{\mathbb{P}}
\newcommand{\XXX}{\mathbb{X}}
\newcommand{\YYY}{\mathbb{Y}}
\newcommand{\Var}{\operatorname{Var}}
\newcommand{\Cov}{\operatorname{Cov}}

\newcommand{\Aa}{\mathcal{A}}
\newcommand{\Bb}{\mathcal{B}}

\newcommand{\eps}{\varepsilon}

\newcommand{\todistr}{\mathop{\longrightarrow}^{d}_{n\to\infty}}
\newcommand{\tofd}{\mathop{\longrightarrow}^{f.d.d.}_{n\to\infty}}
\newcommand{\stofd}{\mathop{\longrightarrow}^{\mathrm{f.d.d.}}_{n\to\infty}}

\newcommand{\stofdn}{\mathop{\longrightarrow}^{\mathrm{f.d.d.}}}

\newcommand{\ome}{t}
\newcommand{\Ome}{T}

\newcommand{\HP}{\Ome}

\newcommand{\ST}{\Ome}
\newcommand{\Ass}{\Ome}

\newcommand{\Pol}{\Ome}

\newtheorem{theorem}{Theorem} 
\newtheorem{lemma}{Lemma} 
\newtheorem{corollary}{Corollary} 
\newtheorem{proposition}{Proposition} 
\newproclaim{remark}{Remark} 

\makeatother

\begin{document}
\begin{frontmatter}

\title{Distribution of levels in high-dimensional random~landscapes}
\runtitle{Distribution of levels in random landscapes}

\begin{aug}
\author{\fnms{Zakhar} \snm{Kabluchko}\corref{}\ead[label=e1]{zakhar.kabluchko@uni-ulm.de}}
\runauthor{Z. Kabluchko}
\affiliation{Ulm University}
\address{Institute of Stochastics\\
Ulm University\\
Helmholtzstr. 18\\
89069 Ulm\\
Germany\\
\printead{e1}}
\end{aug}

\received{\smonth{5} \syear{2010}}
\revised{\smonth{12} \syear{2010}}

%
\begin{abstract}
We prove empirical central limit theorems for the distribution of
levels of various random fields defined on high-dimensional discrete
structures as the dimension of the structure goes to $\infty$. The
random fields considered include costs of assignments, weights of
Hamiltonian cycles and spanning trees, energies of directed polymers,
locations of particles in the branching random walk, as well as
energies in the Sherrington--Kirkpatrick and Edwards--Anderson models.
The distribution of levels in all models listed above is shown to be
essentially the same as in a stationary Gaussian process with regularly
varying nonsummable covariance function. This type of behavior is
different from the Brownian bridge-type limit known for independent or
stationary weakly dependent sequences of random variables.
\end{abstract}

%


\end{frontmatter}

\section{Statement of results}\label{secintro}
\subsection{Introduction}\label{sec1.1}
Strongly correlated random fields defined on high-dimensional discrete
structures arise naturally in stochastic combinatorial optimization and
in the physics of disordered systems.
We will be interested in the properties of the empirical process formed
by the levels of such random fields. The general setting is as follows.
For every $n\in\N$, let $\{\XXX_n(\ome); \ome\in\Ome_n\}$ be a
zero-mean, unit-variance random field with a finite index set $\Ome_n$.
The empirical distribution function of the field $\XXX_n$ counts the
proportion of values of $\XXX_n$ which are not greater than a given
number $z\in\R$. It is defined as
%
\begin{equation}\label{eqdefempiricalfunc}
F_n(z)=\frac{1}{|T_n|}\sum_{\ome\in\Ome_n}1_{\XXX_n(\ome)\leq z},
\qquad z\in\R.
\end{equation}
Here, $|T_n|$ denotes the cardinality of the finite set $T_n$. For a
number of models of stochastic combinatorial optimization we will
prove\vadjust{\goodbreak}
an empirical central limit theorem of the following form:
%
\begin{equation}\label{eqmaineq}
\bigl\{ c_n\bigl(F_n(z)- \E F_n(z)\bigr) ; z\in\R\bigr\}
\stofd
\{p(z)W; z\in\R\}.
\end{equation}
Here, $c_n$ is a normalizing sequence, $\displaystyle\stofdn$ denotes\vspace*{-1pt} the weak
convergence of the finite-dimensional distributions, $p(z)=(2\pi
)^{-1/2}e^{-z^2/2}$ is the standard Gaussian density and $W$ is a
random variable. Both $c_n$ and $W$ depend on the model under
consideration, $W$ being usually normal.



\subsection{Distribution of weights of subgraphs}\label{secsubgraphs}
Our first result deals with the stochastic assignment problem. In this
model, $n$ jobs have to be assigned in a~bijective way to $n$ machines.
The set of all assignments is denoted by $\Ass_n$ and is identified
with the set of all permutations of $n$ elements, so that $|\Ass
_n|=n!$. Let the cost of assigning a job $i$ to the machine $j$ be $\xi
_{i,j}$, where $\{\xi_{i,j}; i,j\in\{1,\ldots,n\}\}$ are independent
copies of a random variable~$\xi$ satisfying $\E\xi=0$ and $\E\xi
^2=1$. The (normalized) cost of an assignment $\ome=(\ome
(i))_{i=1}^n\in\Ass_n$ is then defined by $\XXX_{n}(\ome)=\frac
{1}{\sqrt n}\sum_{i=1}^n \xi_{i, \ome(i)}$. 
%
\begin{theorem}\label{theomainAss}
Let $\{\XXX_n(\ome); \ome\in\Ass_n\}$ be the random landscape of the
stochastic assignment problem. Then, $\XXX_n$ satisfies the empirical
central limit theorem~\eqref{eqmaineq} with $c_n=\sqrt n$ and $W\sim N(0,1)$.
\end{theorem}

The next model we will consider is the mean-field stochastic traveling
salesman problem. Denote by $G_n=(V_n,E_n)$ the undirected complete
graph on a set $V_n$ of $n\geq3$ vertices with the set of edges $E_n$.
A Hamiltonian path in~$G_n$ is a nonoriented closed path which contains
every vertex of~$G_n$ exactly once. Let $\HP_n$ be the set of
Hamiltonian paths in $G_n$. Let the weight of an edge $e\in E_n$ be $\xi
_e$, where $\{\xi_e; e\in E_n\}$ are independent copies of a~random
variable~$\xi$ satisfying $\E\xi=0$ and $\E\xi^2=1$. The (normalized)
weight of a~Hamiltonian path $\ome\in\HP_n$ is then defined by $\XXX
_{n}(\ome)=\frac{1}{\sqrt n}\sum_{e\in\ome} \xi_{e}$.
\begin{theorem}\label{theomainHP}
Let $\{\XXX_n(\ome); \ome\in\HP_n\}$ be the random landscape of the
mean-field stochastic traveling salesman problem. Then, $\XXX_n$
satisfies the empirical central limit theorem~\eqref{eqmaineq} with
$c_n=\sqrt{n/2}$ and $W\sim N(0,1)$.
\end{theorem}

In the next theorem we will deal with the distribution of the weights
of spanning trees in the complete graph. As above, let $G_{n}$ be the
undirected complete graph on~$n$ vertices with the set of edges denoted
by $E_{n}$. A spanning tree is a connected subgraph of $G_{n}$ which
contains all the vertices of $G_{n}$ and has no cycles. Note that the
number of edges in any spanning tree of $G_{n}$ is $n-1$. Let $\ST_n$
be the set of all spanning trees of the complete graph $G_{n}$, the
cardinality of $\ST_n$ being $n^{n-2}$ by the Cayley formula. Let the
weight of an edge $e\in E_{n}$ be $\xi_e$, where $\{\xi_e; e\in E_{n}\}
$ are independent copies of a~random variable~$\xi$ satisfying $\E\xi
=0$ and $\E\xi^2=1$. The (normalized) weight of a spanning tree $\ome
\in\ST_n$ is then defined by $\XXX_{n}(\ome)=\frac{1}{\sqrt{n-1}}\sum
_{e\in\ome} \xi_{e}$.
\begin{theorem}\label{theomainST}
Let the random landscape $\{\XXX_n(\ome); \ome\in\ST_n\}$ representing
the weights of spanning trees be defined as above. Then, $\XXX_n$
satisfies the empirical central limit theorem~\eqref{eqmaineq} with
$c_n=\sqrt{n/2}$ and $W\sim N(0,1)$.
\end{theorem}
%
%
\begin{remark}
We believe that in all our results, the weak convergence of the
finite-dimensional distributions can be replaced by the weak
convergence in the Skorokhod space, but we will not deal with tightness
questions here.
\end{remark}
%
%
\begin{remark}\label{remphin}
In the setting of Theorems~\ref{theomainAss}--\ref{theomainST},
$\lim_{n\to\infty} \E F_n(z)=\Phi(z)$ by the central limit theorem,
where $\Phi(z)$
is the standard Gaussian distribution function. However, we cannot
replace $\E F_n(z)$ by $\Phi(z)$ in~\eqref{eqmaineq}. In order to
justify such a replacement, a relation of the form $\E F_n(z)-\Phi
(z)=o(1/c_n)$ as $n\to\infty$ would be needed. This relation is not
true in general. If the distribution of $\xi$ is nonlattice and $\E|\xi
|^3<+\infty$, then we have, by \cite{gnedenkokolmogorovbook}, page~210,
%
\begin{equation}\label{eqesseen}
\E F_n(z)-\Phi(z)
=
n^{-1/2}Q(z)p(z)+o(n^{-1/2}), \qquad   n\to\infty,
\end{equation}
where $Q(z)=\frac{1}{6}\E[\xi^3](1-z^2)$.
In this case, Theorem~\ref{theomainAss} can be written in the form
%
\begin{equation}\label{eqasymptclt}
\bigl\{\sqrt n\bigl(F_n(z)- \Phi(z)\bigr); z\in\R\bigr\}
\stofd
\bigl\{p(z)\bigl(W+Q(z)\bigr); z\in\R\bigr\},
\end{equation}
where $W\sim N(0,1)$.
Similar considerations apply to Theorems~\ref{theomainHP}, \ref{theomainST}, \ref{theomainBRW}, as well as to the case $d\geq3$
of Theorem~\ref{theomainPol}.
\end{remark}
%


\subsection{Distribution of energies of directed polymers}\label{secpolymers}
A $d$-dimensional directed polymer of length $n$ is a sequence $\ome
=(\ome(k))_{k=0}^{n}$ of sites in $\Z^d$ such that \mbox{$\ome(0)=0$}, and
$\ome(k)$ and $\ome(k+1)$ are neighboring sites for all $k=0,\ldots,
n-1$. The set of all polymers of length $n$ is denoted by $\Pol_n$ and
contains $(2d)^{n}$ elements. Let $\{\xi_{k}(x); k\in\N, x\in\Z^d\}$
be independent copies of a random variable~$\xi$ satisfying $\E\xi=0$
and $\E\xi^2=1$. For $d=1,2$, we additionally assume that $\E|\xi
|^{2+\delta}<\infty$ for some $\delta>0$. The (normalized) energy of a
polymer $\ome\in\Pol_n$ is defined by $\XXX_n(\ome)=\frac{1}{\sqrt
n}\sum_{k=1}^{n} \xi_{k}(\ome(k))$.
\begin{theorem}\label{theomainPol}
Let $\{\XXX_n(\ome); \ome\in\Pol_n\}$ be the random energy landscape of
the directed polymer model. Then, $\XXX_n$ satisfies the empirical
central limit theorem~\eqref{eqmaineq} with
%
\begin{equation}
c_n=
\cases{
\sqrt[4] {\pi n/4}, &\quad$d=1,$\vspace*{2pt}\cr
\sqrt{\pi n/\log n}, &\quad$d=2,$\vspace*{2pt}\cr
\sqrt{n}, &\quad$d\geq3.$
}
\end{equation}
For $d=1,2$, we have $W\sim N(0,1)$. For $d\geq3$, the random variable
$W$ has the same distribution as
$-\sum_{k=1}^{\infty}\sum_{x\in\Z^d} p_{k}(x)\xi_{k}(x),$
where $p_k(x)$ is the probability that a simple $($nearest-neighbor$)$
random walk on $\Z^d$ starting at the origin is at $x\in\Z^d$ at time
$k\in\N$.
\end{theorem}
%
%
\begin{remark}
If $\E|\xi|^3<\infty$, then $|\E F_n(z)-\Phi(z)|\leq C/\sqrt n$ by the
Berry--Esseen inequality (see, e.g., \cite{petrovbook}, page 111).
This implies that we can replace~$\E F_n(z)$ by $\Phi(z)$ in~\eqref{eqmaineq} for $d=1,2$. Note that this does not apply to the case
$d\geq3$. In this case, we may use expansion~\eqref{eqesseen} as in
Remark~\ref{remphin}.
\end{remark}

\subsection{Distribution of particles in the branching random
walk}\label{secBRW}
Branching random walk is a model combining a Galton--Watson branching
process with a random spatial motion of particles. At time $0$ there is
a single particle on the real line located at $0$. At time $1$, this
particle is replaced by a random number of offsprings whose
displacements relative to the position of the parent particle are
i.i.d. random variables. Then, every offspring generates new particles
according to the same rules, and so on. All the random mechanisms
involved are independent.

The formal definition is as follows. Let $\VVV=\bigcup_{n=0}^{\infty}\N^n$
be an infinite tree with root $\varnothing$ (we agree that $\N^0=\{
\varnothing\}$), vertices of the form $t=(v_1,\ldots,v_n)$, where $v_i\in
\N$ and $n=0$ corresponds to the root $t=\varnothing$ and edges
connecting each such~$t$ with its successors $(v_1,\ldots,v_n,k)$,
where $k\in\N$. The number \mbox{$l(t)=n$} is called the length of
$t=(v_1,\ldots,v_n)$. Let $\{Z_t;t\in\VVV\}$ be independent copies of
a random variable $Z$ which takes values in $\N$ and satisfies $m:=\E
Z>1$ and $\E Z^2<\infty$. The random variable $Z_t$ should be thought
of as the number of children of the particle coded by the vertex $t$.
The $n$th generation of the branching random walk is the random set
$T_n$ consisting of all vertices $t=(v_1,\ldots,v_n)$ of length $n\in\N
$ such that $v_{k}\leq Z_{(v_1,\ldots,v_{k-1})}$ for every $k=1,\ldots,n$.
Independently of the $Z_t$'s, let $\{\xi_t; t\in\VVV\setminus\{
\varnothing\}\}$ be independent copies of a random variable $\xi$ such
that $\E\xi=0$, $\E\xi^2=1$. The random variable~$\xi_t$ should be
thought of as the displacement of the particle coded by the vertex $t$
relative to its parent. For $t=(v_1,\ldots,v_n)\in\VVV\setminus\{
\varnothing\}$ define $\XXX_n(t)=\frac{1}{\sqrt{n}}\sum_{k=1}^n \xi
_{(v_1,\ldots,v_{k})}$. Then, $\{\XXX_n(\ome); \ome\in\Ome_n\}$ are the
normalized positions of the particles in the $n$th generation of the
branching random walk.

\begin{theorem}\label{theomainBRW}
The random field $\{\XXX_n(t), t\in T_n\}$ defined as above satisfies
the empirical central limit theorem~\eqref{eqmaineq} with $c_n=\sqrt
n$. The limiting random variable~$W$ has the same distribution as
$-\!\lim_{n\to\infty}\sqrt n |T_n|^{-1}\sum_{t\in T_n} \XXX_n(t)
$.
\end{theorem}
%

In the case of Bernoulli-distributed displacements this theorem is due
to~\cite{chen01}.
The method of~\cite{chen01} relies strongly on the Markov property of
the branching random walk. We will recover Theorem~\ref{theomainBRW}
as a particular case of our general approach.

\subsection{Distribution of energy levels in spin glasses}\label{secgauss}
Our last result concerns the distribution of energy levels in spin
glasses. The general setting is as follows.
For every $n\in\N$, let $G_n=(V_n,E_n)$ be an undirected graph without
loops and multiple edges on a finite set of vertices $V_n$ with the set
of edges $E_n$. A spin configuration is a map $t\dvtx V_n\to\{-1,1\}$. Let
$\Ome_n=\{-1,1\}^{V_n}$ be the set of all spin configurations. Spins
located at vertices $v_1$ and $v_2$ interact if there is an edge $e=\{
v_1,v_2\}\in E_n$, the energy of the interaction being $\ome(v_1)\ome
(v_2) J(e)$, where $\{J(e);e\in E_n\}$ are independent standard
Gaussian random variables. The energy of a spin configuration $\ome\in
\Ome_n$ is defined by
%
\begin{equation}\label{eqdefEA}
\XXX_n(\ome)=|E_n|^{-1/2}\sum_{e=\{v_1,v_2\}\in E_n} \ome(v_1)\ome(v_2) J(e).
\end{equation}
Examples are provided by the Sherrington--Kirkpatrick model in which
$G_n$ is the complete graph on $n$ vertices, and the $d$-dimensional
Edwards--Anderson model, in which $G_n$ is the $d$-dimensional discrete
box with side length $n$ and nearest-neighbor interactions.
\begin{theorem}\label{theomainEA}
Let $\{\XXX_n(\ome); \ome\in\Ome_n\}$ be the energy landscape defined
as in~\eqref{eqdefEA}. If $\lim_{n\to\infty}|E_n|=\infty$, then the
following empirical central limit theorem holds:
%
\begin{equation}
\biggl\{|E_n|^{1/2}2^{-|V_n|}\sum_{\ome\in\Ome_n}\bigl(1_{\XXX_n(\ome)\leq
z}-\Phi(z)\bigr); z\in\R\biggr\}
\tofd
\biggl\{\frac{zp(z)}{\sqrt2}W; z\in\R\biggr\},
\end{equation}
where $\Phi$ is the standard Gaussian distribution function and $W\sim N(0,1)$.
\end{theorem}
%


\subsection{Discussion}\label{sec1.6}
Empirical central limit theorems have been extensively studied for
stationary sequences of random variables under various short-range
dependence conditions. For example, it has been shown in~\cite{sun65,berman70} that if $\{\XXX(n); n\in\Z\}$ is a stationary
zero-mean, unit-variance Gaussian process whose covariance function
$r(n)=\E[\XXX(0)\XXX(n)]$ satisfies $\sum_{n\in\Z}|r(n)|<\infty$, then
%
\begin{equation}
\Biggl\{\frac{1}{\sqrt n} \sum_{k=1}^n \bigl(1_{\XXX(k)\leq z}-\Phi(z)\bigr); z\in
\R\Biggr\}\stofd\{B(z); z\in\R\},
\end{equation}
where
$\{B(z); z\in\R\}$ is a zero-mean Gaussian process with covariance function
%
\begin{equation}\label{eqcorrberman}
\Cov(B(z_1),B(z_2))
=
\sum_{k\in\Z}\Cov\bigl(1_{\XXX(0)\leq z_1}, 1_{\XXX(k)\leq z_2}\bigr),
\qquad z_1,z_2\in\R.
\end{equation}
Similar results are available for stationary processes under mixing
conditions~\cite{billingsleybook}, Chapter 22, \cite{deo73,withers75},
stationary associated sequences~\cite{yu93}, to cite only a few
references. 

There has been also much interest in proving empirical central limit
theorems for stationary long-range dependent processes (see~\cite{taqqu74,dobrushinmajor79,taqqu79,dehlingtaqqu}, as well as the
monographs~\cite{ivanovleonenkobook,leonenkobook,dehlingbook} for
further references).
It has been shown that if $\{\XXX(n);  n\in\Z\}$ is a stationary
zero-mean, unit-variance Gaussian process whose covariance function $r$
satisfies $r(n)=L(n) n^{-D}$ for some function~$L$ that varies slowly
at $+\infty$ and some $D\in(0,1)$, then
%
\begin{equation}\label{eqCLTGaussRV}
\Biggl\{\frac{C_D}{L^{1/2}(n) n^{1-(1/2) D}} \sum_{k=1}^n \bigl(1_{\XXX
(k)\leq z}-\Phi(z)\bigr); z\in\R\Biggr\}
\stofd
\{p(z)W; z\in\R\},\hspace*{-30pt}
\end{equation}
where $C_D>0$ is some explicit constant and $W\sim N(0,1)$ (see~\cite{taqqu74,dehlingtaqqu} for stronger results).

The models considered in the present paper look, at a first sight,
rather different from stationary Gaussian processes with regularly
varying covariance function. Nevertheless, as far as the empirical
process is concerned, they behave in essentially the same way as
in~\eqref{eqCLTGaussRV}. A nonrigorous explanation of this
phenomenon will be given in Section~\ref{secheuristic}. 

Let us also mention that several authors proved Poisson limit theorems
for the local distribution of values of highly-correlated random fields
in small windows~\cite{baukemertens,bovierkurkova07,benarouskuptsov09,borgsetal1,borgsetal2}.
As opposed to these results, we consider the distribution of values of
random fields on a global scale.

%

\subsection{Idea of the proofs}\label{secheuristic}
Let us describe a nonrigorous argument justifying our results. As an
approximation to the models considered in \mbox{Theorems~\ref{theomainAss}--\ref{theomainBRW}}, we take $\{\XXX_n(\ome); \ome\in\Ome_n\}$ to be a
Gaussian process with zero-mean, unit-variance marginals and a covariance
structure given by $\E[\XXX_n(\ome_1)\XXX_n(\ome_2)]=\eps_n$ for all
$\ome_1\neq\ome_2$, where $\eps_n\in(0,1)$ is some sequence tending
to $0$ as $n\to\infty$. Intuitively, the sequence $\eps_n$ represents
the order of the overlap of two generic assignments, Hamiltonian paths, etc.

Let $F_n$ be the empirical distribution function of $\XXX_n$ defined as
in~\eqref{eqdefempiricalfunc}.
The process~$\XXX_n$ can be represented (in distribution) as $\XXX
_{n}(\ome)=\sqrt{1-\eps_n} \XXX'_{n}(\ome)+\sqrt{\eps_n} N$,
where $\{\XXX_n'(\ome); \ome\in\Ome_n\}$ and $N$ are independent
standard Gaussian random variables. Thus, we have a representation
\[
F_n(z)=F_n'\biggl(\frac{z-\sqrt{\eps_n} N}{\sqrt{1-\eps_n}}\biggr),\qquad
 z\in\R, 
\]
where $F'_n(z)=
\frac1{|\Ome_n|} \sum_{\ome\in\Ome_n} 1_{\XXX_n'(\ome)\leq z}$ is the
empirical distribution function of~$\XXX_n'$.
By the central limit theorem, we have $F_n'\approx\Phi$ as $n\to\infty
$ with a Brownian bridge error term of order $1/\sqrt{|\Ome_n|}$,
where $\Phi$ is the standard Gaussian distribution function. Now, the
common feature of the models considered in Theorems~\ref{theomainAss}--\ref{theomainBRW} is that $\eps_n$, the order of the
correlation of two generic elements in $\Ome_n$, is much larger than
$1/|\Ome_n|$. So, the order of the Brownian bridge fluctuations is much
smaller than the order of the shift $\sqrt{\eps_n}  N$. Thus, we may write
%
\begin{equation}
F_n(z)= F_n'\biggl(\frac{z-\sqrt{\eps_n} N}{\sqrt{1-\eps_n}}\biggr)
\approx
\Phi\biggl(\frac{z-\sqrt{\eps_n} N}{\sqrt{1-\eps_n}}\biggr)
\approx
\Phi(z)-\sqrt{\eps_n} p(z)N,\hspace*{-25pt}
\end{equation}
where $\approx$ means that we are ignoring terms of order $o_P(\sqrt
{\eps_n})$ as $n\to\infty$.
This leads to a result of the form
%
\begin{equation}
\biggl\{\frac1 {\sqrt{\eps_n}} \bigl(F_n(z)-\Phi(z)\bigr); z\in\R\biggr\} \stofd
\{-p(z)N; z\in\R\}.
\end{equation}

For example, let $\{\XXX(n); n\in\Z\}$ be a stationary zero-mean,
unit-variance Gaussian process whose covariance function $r$ satisfies
$r(n)=L(n)n^{-D}$, where $L$ is a slowly varying function and $D>0$.
Then, for generic $k_1,k_2$ in $\Ome_n=\{1,\ldots,n\}$, $\Cov(\XXX
(k_1),\XXX(k_2))$ is of order $\eps_n\approx L(n)n^{-D}$. If $D\in
(0,1)$, then $\eps_n$ is asymptotically larger than $1/|\Ome_n|$ and
the heuristic applies; cf.~\eqref{eqCLTGaussRV}. In the models of
Section~\ref{secsubgraphs} and in the branching random walk, we have
$\eps_n\approx1/n$, whereas $|\Ome_n|$ grows exponentially, so that
again $\eps_n$ is larger than $1/|\Ome_n|$. For directed polymers, $\eps
_n$ depends on the dimension $d$ and is again larger than $1/|\Ome_n|$.


On a more rigorous level, our proofs will be based on an adaptation of
the reduction method of \cite{taqqu74}. This method was introduced in
the setting of stationary Gaussian processes with regularly varying
covariance. The idea is to approximate the empirical distribution
function by a certain expansion involving Hermite polynomials. Recall
that the Hermite polynomials form an orthogonal system with respect to
the weight $p$, the standard Gaussian density; see Section~\ref{sechermitepoly} for precise definitions. Every function which is
square integrable with respect to the weight $p$ can be expanded into a
Hermite--Fourier series. For the function $f(x)=1_{x\leq z}-\Phi(z)$
(here, $z\in\R$ is fixed), the first two terms in the Hermite--Fourier
expansion are
%
\begin{equation}\label{eqfourier-hermite}
f(x)=1_{x\leq z}-\Phi(z)=-p(z)x-\tfrac12 zp(z)(x^2-1)+\cdots.
\end{equation}
To prove Theorems~\ref{theomainAss}--\ref{theomainBRW}, we will show
that the random variable\break $\sum_{t\in T_n} (1_{\XXX_n(\ome)\leq z}-\P
[\XXX_n(\ome)\leq z])$ can be approximated in the $L^2$-sense by the
random variable $-p(z)\sum_{t\in T_n}\XXX_n(\ome)$ corresponding to the
first term of the expansion~\eqref{eqfourier-hermite}. The statements
justifying this approximation are Lemma~\ref{lemest1} and
Proposition~\ref{propmain} below.
For the proof of Theorem~\ref{theomainEA}, we need a more accurate
approximation involving the second Hermite polynomial since there, we
have $\sum_{t\in T_n}\XXX_n(\ome)=0$ by symmetry reasons. In the
setting of Theorem~\ref{theomainEA}, we will prove that $\sum_{t\in
T_n} (1_{\XXX_n(\ome)\leq z}- \Phi(z))$ can be approximated by $-\frac
12 zp(z)\sum_{t\in T_n}(\XXX_n^2(\ome)-1)$.





\subsection{Notation}\label{sec1.8}
Let us collect the notation which will be used throughout the paper.
The standard Gaussian density and distribution function are denoted by\vadjust{\goodbreak}
$p(z)=(2\pi)^{-1/2}e^{-z^2/2}$ and $\Phi(z)=\int_{-\infty}^{z}p(t)\,dt$,
respectively.
We denote by $\xi$ a random variable satisfying $\E\xi=0$ and $\E\xi
^2=1$. Let $\Phi_n$ be the distribution function of $(\xi_1+\cdots+\xi
_n)/\sqrt n$, where $\{\xi_i; i\in\N\}$ are independent copies of $\xi
$. By the central limit theorem, $\lim_{n\to\infty}\Phi_n(z)=\Phi(z)$
for every $z\in\R$. Throughout, $C$ is a large positive constant whose
value may change from line to line.

\section{Proofs for combinatorial models}\label{secproof}
\subsection{Local limit theorems}\label{sec2.1}
We start by recalling two classical local limit theorems which will be
needed in our proofs. The first of them deals with lattice random
variables. Recall that a random variable is called lattice if its
values are of the form $b+h\Z$ for some $b\in\R$ and $h\geq0$.
\begin{theorem}[(\cite{gnedenkokolmogorovbook}, page~233, or~\cite{petrovbook}, page~187)]\label{theognedenko}
Let $\{\xi_i; i\in\N\}$ be independent copies of a random variable $\xi
$ satisfying $\E\xi=0$ and $\E\xi^2=1$. Assume that the values of $\xi
$ are of the form $b+h\Z$, where $h>0$ is maximal with this property.
Then, the following asymptotic relation holds uniformly in $z\in nb+h\Z$:
%
\begin{equation}
\P[\xi_1+\cdots+\xi_n=z]
=
\frac{h}{\sqrt n}p\biggl(\frac{z}{\sqrt n}\biggr)+o\biggl(\frac{1}{\sqrt
n}\biggr),
 \qquad   n\to\infty.
\end{equation}
\end{theorem}

The next theorem is an analogue of Theorem~\ref{theognedenko} for
nonlattice distributions. Recall the notation introduced at the end of
Section~\ref{secintro}.
\begin{theorem}[(\cite{stone65})]\label{theostone}
Let $\xi$ be a nonlattice random variable satisfying\break $\E\xi=0$ and $\E
\xi^2=1$. Then, the following asymptotic relation holds uniformly in
$z_1, z_2\in\R$:
%
\begin{eqnarray}\label{eqthmstone}
\Phi_n(z_2)-\Phi_n(z_1)=\Phi(z_2)-\Phi(z_1)+o(1)(|z_2-z_1|+n^{-1/2}),\\
  \eqntext{ n\to\infty.}
\end{eqnarray}
\end{theorem}
\begin{corollary}\label{corstone}
Regardless of whether $\xi$ is lattice or nonlattice, there is a~constant $C>0$ depending on $\xi$ such that for all $n\in\N$ and
$z_1,z_2\in\R$,
%
\begin{equation}
|\Phi_n(z_2)-\Phi_n(z_1)|\leq C|z_2-z_1|+Cn^{-1/2}.
\end{equation}
\end{corollary}
\begin{pf} 
If the distribution of $\xi$ is nonlattice, then the corollary follows
immediately from~\eqref{eqthmstone} and the fact that the function
$\Phi$ is Lipschitz. Suppose that $\xi$ is lattice as in Theorem~\ref{theognedenko}. Without restriction of generality, let $z_1<z_2$ and
define $I_n=(nb+h\Z) \cap(\sqrt n z_1,\sqrt n z_2]$. Then by
Theorem~\ref{theognedenko},
\[
\Phi_n(z_2)-\Phi_n(z_1)
=
\sum_{z\in I_n} \biggl(\frac{h}{\sqrt n}p\biggl(\frac{z}{\sqrt n}
\biggr)+o\biggl(\frac{1}{\sqrt n}\biggr)\biggr)
\leq
\sum_{z\in I_n} \frac C{\sqrt n},
\]
where the $o$-term is uniform in $z\in\R$.
Since the cardinality of $I_n$ differs from $h^{-1}\sqrt n (z_2-z_1)$
by at most $1$, we obtain the statement of the corollary.
\end{pf}
\begin{remark}
With an additional assumption $\E|\xi|^3<\infty$, Corollary~\ref{corstone} follows from the Berry--Esseen inequality (see~\cite{petrovbook}, page 111).
\end{remark}

\subsection{The main lemma}\label{sec2.2}
The next lemma will play a crucial role in the sequel. Essentially, it
provides an estimate for the dependence between the random variables
$1_{X_1\leq z}$ and $1_{X_2\leq z}$, where $X_1$ and $X_2$ are two
normalized sums of i.i.d. random variables having a nontrivial
overlap. In our applications, $X_1$ and $X_2$ will be the normalized
weights of two Hamiltonian paths, spanning trees, etc. We will
regularize $1_{X_1\leq z}$ and $1_{X_2\leq z}$ by subtracting certain
terms motivated by the Hermite expansion of the function $f(x)=1_{x\leq z}$.

\begin{lemma}\label{lemest1}
Let $\{\xi_{i}; i\in\N\}$ be independent copies of a random variable
$\xi$ satisfying $\E\xi=0$ and $\E\xi^2=1$. Let $z\in\R$ be fixed.
Given $r\in\N\cup\{0\}$, $n\in\N$ with $r\leq n$, define two random
variables $Y_1,Y_2$ by
%
\begin{equation}\label{eqdefXin}
Y_i=1_{X_i\leq z}- \Phi_n(z)+p(z)X_i, \qquad   i=1,2,
\end{equation}
where $X_{1}=\frac{1}{\sqrt n}\sum_{i=1}^n \xi_i$ and $X_{2}=\frac
{1}{\sqrt n}\sum_{i=n-r+1}^{2n-r} \xi_i$.
Then, there is a constant~$C$ depending only on the distribution of $\xi
$ such that for all $r\in\N\cup\{0\}$, $n\in\N$ with $r\leq n$, we
have
%
\begin{equation}\label{eqlemp1}
0\leq\E[Y_1Y_2]\leq C\frac{r}{n}.
\end{equation}
Further, if $\eps_n>0$ is any sequence with $\lim_{n\to\infty} \eps
_n=0$, then there is a sequence $\delta_n$ such that $\lim_{n\to\infty
}\delta_n=0$ and for every $r\in\N\cup\{0\}$, $n\in\N$ with $r\leq\eps
_n n$, we have
%
\begin{equation}\label{eqlemp2}
0\leq\E[Y_1Y_2]\leq\delta_n\frac r n.
\end{equation}
\end{lemma}
\begin{pf}
Since the statement is trivially fulfilled for $r=0$ and $r=n$, we
assume $0<r<n$ henceforth. It will be convenient to introduce the
following notation: for $u\in\R$, we write
%
\begin{equation}\label{eqdefrho}
\rho=\frac r n\in(0,1),\qquad
z(u)=\frac{z-u\sqrt{\rho}}{\sqrt{1-\rho}}.
\end{equation}
It follows from~\eqref{eqdefXin} that we have
%
\begin{equation}\label{eqcovXij}
\E[Y_1Y_2]
=\Cov(1_{X_1\leq z}, 1_{X_2\leq z})+2p(z)\E[1_{X_1\leq z}
X_2]+p^2(z)\rho.
\end{equation}
We start by considering the first term on the right-hand side of~\eqref{eqcovXij}. We are going to show that
\begin{eqnarray}\label{eqcovdelta1}
&&\Cov(1_{X_1\leq z}, 1_{X_2\leq z})
\nonumber
\\[-8pt]
\\[-8pt]
\nonumber
&&\qquad=
\frac12 \int_{\R}\int_{\R}\bigl(\Phi_{n-r}(z(u))-
\Phi_{n-r}(z(v))
\bigr)^2\Phi_r(du)\Phi_r(dv).
\end{eqnarray}
Define three independent random variables $\tilde X_1$, $\tilde X$,
$\tilde X_2$ by
%
\begin{equation}
\tilde X_1=\frac1 {\sqrt n}\sum_{i=1}^{n-r}\xi_i,\qquad
\tilde X=\frac1 {\sqrt n}\sum_{i=n-r+1}^n\xi_i,\qquad
\tilde X_2=\frac1 {\sqrt n}\sum_{i=n+1}^{2n-r}\xi_i.
\end{equation}
Note that $X_1=\tilde X_1+\tilde X$ and $X_2=\tilde X_2+\tilde X$.
The distribution function of $\tilde X/\sqrt{\rho}$ is $\Phi_r$.
Conditioning on the event $\tilde X/\sqrt{\rho}\in du$ and using the
independence of $\tilde X_1,\tilde X,\tilde X_2$, we obtain
\begin{eqnarray}\label{eqlemwsp1}
\E[1_{X_1\leq z}1_{X_2\leq z}]
&=&
\P[\tilde X_1+\tilde X\leq z, \tilde X_2+\tilde X\leq z
]\nonumber\\
&=&
\int_{\R}\bigl(\P\bigl[\tilde X_1\leq z-u\sqrt{\rho} \bigr]\bigr)^2
\Phi_{r}(du)
\nonumber
\\[-8pt]
\\[-8pt]
\nonumber
&=&
\int_{\R} \Phi_{n-r}^2(z(u))\ \Phi_{r}(du)\\
&=&
\frac12
\int_{\R}\int_{\R}
\bigl(
\Phi_{n-r}^2(z(u))+\Phi_{n-r}^2(z(v))
\bigr)
\Phi_r(du)\Phi_r(dv).\nonumber
\end{eqnarray}
In a similar way, we obtain
\begin{eqnarray}\label{eqlemwsp2}
\E[1_{X_1\leq z}]\E[1_{X_2\leq z}]
&=&
(\P[\tilde X_1+\tilde X\leq z ])^2\nonumber\\
&=&
\biggl(\int_{\R}\Phi_{n-r}(z(u))\Phi_r(du)\biggr)^2\\
&=&
\int_{\R}\int_{\R}
\Phi_{n-r}(z(u))
\Phi_{n-r}(z(v))
\Phi_r(du)\Phi_r(dv).\nonumber
\end{eqnarray}
Bringing~\eqref{eqlemwsp1} and~\eqref{eqlemwsp2} together, we
obtain~\eqref{eqcovdelta1}.
Let us consider the second term on the right-hand side of~\eqref{eqcovXij}. Conditioning on $\tilde X/\sqrt{\rho}\in du$, we obtain
\begin{eqnarray}\label{eqwspom5}
&&2p(z)\E[1_{X_1\leq z} X_2]\nonumber\\
&&\qquad=
2p(z)\E[1_{\tilde X_1+\tilde X\leq z} \tilde X]\nonumber\\
&&\qquad=
2p(z)\sqrt{\rho} \int_{\R} u \P\bigl[\tilde X_1\leq z-u\sqrt{\rho}
\bigr]\Phi_r(du)\\
&&\qquad=
2p(z)\sqrt{\rho} \int_{\R}u\Phi_{n-r}(z(u))\Phi_r(du)\nonumber\\
&&\qquad=p(z)\sqrt{\rho} \int_{\R} \int_{\R} (u-v) \bigl(\Phi_{n-r}(z(u)
)- \Phi_{n-r}(z(v))\bigr)\Phi_r(du)\Phi_r(dv).\nonumber
\end{eqnarray}
Also, we have
%
\begin{equation}\label{eqwspom4}
\frac12 \int_{\R} \int_{\R} (u-v)^2\Phi_r(du)\Phi_r(dv)=1.
\end{equation}
Bringing~\eqref{eqcovXij}, \eqref{eqcovdelta1}, \eqref{eqwspom5},
\eqref{eqwspom4} together, we obtain
%
\begin{equation}\label{eqcovXdelta}
\E[Y_1Y_2]
=
\frac12 \int_{\R}\int_{\R}\Delta^2(u,v)
\Phi_r(du)\Phi_r(dv),
\end{equation}
where $\Delta(u,v)$ is given by
%
\begin{equation}\label{eqdefdelta}
\Delta(u,v)
=
\Phi_{n-r}(z(u))
-
\Phi_{n-r}(z(v))
+p(z)(u-v)\sqrt{\rho}.
\end{equation}

Let us now prove the first statement of the lemma. Note that~\eqref{eqcovXdelta}
implies that $\E[Y_1Y_2]\geq0$. It follows from~\eqref{eqdefXin} that $\E Y_1^2 =\E Y_2^2\leq9$. By the Cauchy--Schwarz
inequality, equation~\eqref{eqlemp1} is fulfilled for $r\in[n/2,n]$
and every $n\in\N$ with $C=18$. Let us henceforth assume that $r\leq
n/2$ (and so, $\rho\leq1/2$).
Applying Corollary~\ref{corstone} and recalling~\eqref{eqdefrho}, we
obtain, both in the lattice and in the nonlattice case,
%
\begin{eqnarray}\label{eqlemwsp3}
|\Phi_{n-r}(z(u))-\Phi_{n-r}(z(v))|
&\leq&
C\biggl(\frac{|u-v|\sqrt{\rho}}{\sqrt{1-\rho}}+\frac{1}{\sqrt
{n-r}}\biggr)
\nonumber
\\[-8pt]
\\[-8pt]
\nonumber
&\leq&
C (|u-v|+1)\sqrt{\rho}.
\end{eqnarray}
It follows from~\eqref{eqdefdelta} and~\eqref{eqlemwsp3} that
$|\Delta(u,v)|
\leq C (|u-v|+1) \sqrt{\rho}$.
Hence,
%
\begin{equation}\label{eqintegrDelta}
\Delta^2(u,v)
\leq
C\bigl((u-v)^2+1\bigr)\rho.
\end{equation}
Inserting this into~\eqref{eqcovXdelta} yields
%
\begin{equation}\label{eqlemwsp4}
\E[Y_1Y_2]
\leq
C \rho\int_{\R}\int_{\R} \bigl((u-v)^2+1\bigr)
\Phi_r(du)\Phi_r(dv)=3C\rho.
\end{equation}
This completes the proof of~\eqref{eqlemp1}.

Let us prove the second statement of the lemma. It suffices to show
that for every $\delta>0$ there is $N=N(\delta)$ such that for every
$n>N$ and $r\leq\eps_n n$, we have $\E[Y_1Y_2]\leq\delta\rho$.
It follows from~\eqref{eqintegrDelta}, \eqref{eqwspom4} and the weak
convergence of $\Phi_r$ to $\Phi$ as $r\to\infty$ that we can choose
$B=B(\delta)$ such that for all $n,r\in\N$ with $r\leq n/2$,
%
\begin{equation}\label{eqlemwsp5}
\int_{\R^2\setminus[-B,B]^2}\Delta^2(u,v) \Phi_r(du)\Phi_r(dv)<\delta
\rho.
\end{equation}

Assume first that the distribution of $\xi$ is nonlattice. We always
assume that $r\leq\eps_n n$. By Theorem~\ref{theostone}, the
following holds uniformly in $u,v\in[-B,B]$ as $n\to\infty$:
\begin{eqnarray*}
\Phi_{n-r}(z(u))
-\Phi_{n-r}(z(v))
&=&
p(z)\bigl(z(u)-z(v)\bigr)+o\bigl(\sqrt{\rho}\bigr)\\
&=&
-p(z)\bigl((u-v)+o(1)\bigr)\sqrt{\rho}.
\end{eqnarray*}
Together with~\eqref{eqdefdelta}, this implies that $\Delta
(u,v)=o(\sqrt{\rho})$ uniformly in $u,v\in[-B,B]$ as $n\to\infty$.
It follows that for $n$ large enough,
%
\begin{equation}\label{eqlemwsp6}
\int_{[-B,B]^2}\Delta^2(u,v) \Phi_r(du)\Phi_r(dv)<\delta\rho.
\end{equation}
This, together with~\eqref{eqlemwsp5} and~\eqref{eqcovXdelta},
completes the proof in the nonlattice case.\vadjust{\goodbreak}

Assume now that the random variable $\xi$ is lattice with values in the
set $b+h\Z$, with $h$ being maximal with this property. Let $u,v\in
[-B,B]\cap r^{-1/2}(rb+h\Z)$ with $u<v$. Note that by~\eqref{eqdefrho}, $z(u)-z(v)\in(n-r)^{-1/2}h\Z$. Hence, the number of
points in the set
\[
I_{n,r}(u,v):=(z(v), z(u)]\cap(n-r)^{-1/2}\bigl((n-r)b+h\Z\bigr)
\]
is equal to $h^{-1}(n-r)^{1/2}(z(u)-z(v))$. By Theorem~\ref{theognedenko},
\begin{eqnarray*}
\Phi_{n-r}(z(u))
-\Phi_{n-r}(z(v))
&=&\sum_{x\in I_{n,r}(u,v)} \P\biggl[\frac{\xi_1+\cdots+\xi_{n-r}}{\sqrt
{n-r}}= x\biggr]\\
&=&\sum_{x\in I_{n,r}(u,v)} \biggl(\frac{h}{\sqrt{n-r}}p(x)+o\biggl(\frac
{1}{\sqrt{n-r}}\biggr) \biggr)\\
&=&(v-u)p(z)\sqrt{\rho}+o\bigl(\sqrt{\rho}\bigr).
\end{eqnarray*}
It follows that $\Delta(u,v)=o(\sqrt{\rho})$ as $n\to\infty$ uniformly
in $u,v\in[-B,B]\cap r^{-1/2}(rb+h\Z)$. Hence, equation~\eqref{eqlemwsp6} holds for $n$ large enough and the proof is complete.
\end{pf}

\subsection{An empirical central limit theorem for overlapping sums}\label{sec2.3}
In this section, we state and prove a result from which we will deduce
Theorems~\ref{theomainAss}--\ref{theomainBRW}. It is an empirical central
limit theorem for overlapping sums of independent random variables. Let
$\{\xi_{e}; e\in E\}$ be independent copies of a random variable $\xi$
satisfying $\E\xi=0$ and $\E\xi^2=1$, where $E$ is some countable
index set. For every $n\in\N$, let $\Ome_n\subset2^E$ be a finite
collection of (typically, overlapping) subsets of $E$, each subset
having cardinality $n$. Define a random field $\{\XXX_n(\ome); \ome\in
\Ome_n\}$ by
%
\begin{equation}\label{eqpropdefXn}
\XXX_{n}(\ome)=\frac1 {\sqrt n}\sum_{e\in\ome}\xi_{e}.
\end{equation}
Let $\Phi_n(z)=\P[\XXX_n(t)\leq z]$, where $z\in\R$, be the
distribution function of $\XXX_n(t)$. The covariance function of the
random field $\XXX_n$ is given by $\rho_{n}(\ome_1,\ome_2)=\frac1n
|\ome_1\cap\ome_2|$.
Define also $s_n\geq0$ by
%
\begin{equation}\label{eqpropdef}
s_n^2
=\Var\biggl[\sum_{\ome\in\Ome_n} \XXX_n(\ome)\biggr]
=\sum_{\ome_1,\ome_2\in\Ome_n}\rho_{n}(\ome_1,\ome_2).
\end{equation}

\begin{proposition}\label{propmain}
Let the random field $\{\XXX_n(\ome); \ome\in\Ome_n\}$ be defined as
above. Assume that for some random variable $V$ and some sequence $\eps
_n>0$ with $\lim_{n\to\infty}\eps_n=0$, the following two conditions
are satisfied:
%
\begin{eqnarray}\label{eqpropcond1}
\frac{1}{s_n} \sum_{\ome\in\Ome_n} \XXX_n(\ome)&\displaystyle\mathop{\longrightarrow}^{d}_{n\to\infty}& V,
\\
\label{eqpropcond2}
\lim_{n\to\infty}\frac{1}{s_n^2}
\sum_{\ome_1,\ome_2\in\Ome_n}\rho_{n}(\ome_1,\ome_2)1_{\rho_n(\ome
_1,\ome_2)>\eps_n }\hspace*{-6pt}&=&\hspace*{-6pt}0.
\end{eqnarray}
Then, the following convergence of stochastic processes holds true:
%
\begin{equation}\label{propstatement}
\biggl\{\frac{1}{s_n}\sum_{\ome\in\Ome_n} \bigl(1_{\XXX_n(\ome)\leq z}- \Phi
_n(z)\bigr); z\in\R\biggr\}
\tofd
\{-p(z)V; z\in\R\}.
\end{equation}
\end{proposition}
\begin{pf}
For $z\in\R$, define a zero-mean random field $\{\YYY_n(\ome;z); \ome\in
\Ome_n\}$ by
%
\begin{equation}\label{eqdefYnprop}
\YYY_n(\ome;z)=1_{\XXX_n(\ome)\leq z}- \Phi_n(z)+p(z)\XXX_n(\ome).
\end{equation}
We will show that
%
\begin{equation}\label{eqneedL2}
\quad\lim_{n\to\infty}\Var\biggl[\frac{1}{s_n} \sum_{\ome\in\Ome_n} \YYY
_n(\ome;z)\biggr]
=\lim_{n\to\infty} \frac{1}{s_n^2} \sum_{\ome_1,\ome_2\in\Ome_n} \E[\YYY
_n(\ome_1;z)\YYY_n(\ome_2;z)]
=0.
\end{equation}
By the first part of Lemma~\ref{lemest1}, we have for every $\ome
_1,\ome_2\in\Ome_n$,
%
\begin{equation}\label{eqwspom1}
0\leq\E[\YYY_{n}(\ome_1;z)\YYY_{n}(\ome_2;z)]\leq C\rho_{n}(\ome_1,\ome_2).
\end{equation}
This allows us to estimate the contribution of those terms in~\eqref{eqneedL2} which satisfy $\rho_{n}(\ome_1,\ome_2)>\eps_n$. It follows
from~\eqref{eqwspom1} and~\eqref{eqpropcond2} that as $n\to\infty$,
%
\begin{equation}\label{eqwspom2}
\mathop{\sum_{\ome_1,\ome_2\in\Ome_n}}_{\rho_{n}(\ome_1,\ome_2)>\eps
_n} \E[\YYY_{n}(\ome_1;z)\YYY_{n}(\ome_2;z)]
\leq
C\mathop{\sum_{\ome_1,\ome_2\in\Ome_n}}_{\rho_{n}(\ome_1,\ome_2)>\eps_n
}\rho_{n}(\ome_1,\ome_2)=o(s_n^2).
\end{equation}
Let us consider the terms with $\rho_{n}(\ome_1,\ome_2)\leq\eps_n$. It
follows from the second part of Lemma~\ref{lemest1} that there is a
sequence $\delta_n>0$ such that $\lim_{n\to\infty}\delta_n=0$ and for
every~$\ome_1,\ome_2$ such that $\rho_{n}(\ome_1,\ome_2)\leq\eps_n$,
we have
%
\begin{equation}\label{eqpropprwsp1}
0\leq\E[\YYY_{n}(\ome_1;z)\YYY_{n}(\ome_2;z)]\leq\delta_n \rho
_{n}(\ome_1,\ome_2).
\end{equation}
It follows from~\eqref{eqpropprwsp1} and~\eqref{eqpropdef} that as
$n\to\infty$,
%
\begin{equation}\label{eqwspom3}
\mathop{\sum_{\ome_1,\ome_2\in\Ome_n}}_{\rho_{n}(\ome_1,\ome_2)\leq\eps
_n} \E[\YYY_{n}(\ome_1;z)\YYY_{n}(\ome_2;z)]
\leq
\delta_n\mathop{\sum_{\ome_1,\ome_2\in\Ome_n}}_{\rho_{n}(\ome_1,\ome
_2)\leq\eps_n} \rho_{n}(\ome_1,\ome_2) =o(s_n^2).
\end{equation}
Combining~\eqref{eqwspom2} and~\eqref{eqwspom3}, we obtain~\eqref{eqneedL2}.

Take some $z_1,\ldots,z_d\in\R$. Recalling~\eqref{eqdefYnprop}, we
may write for every $i=1,\ldots,d$,
\[
\frac{1}{s_n}\sum_{\ome\in\Ome_n} \bigl(1_{\XXX_n(\ome)\leq z_i}- \Phi_n(z_i)\bigr)
=
-\frac{p(z_i)}{s_n} \sum_{\ome\in\Ome_n}\XXX_n(\ome)+
\frac{1}{s_n} \sum_{\ome\in\Ome_n} \YYY_n(\ome;z_i).
\]
The first term on the right-hand side converges to $-p(z_i)V$ in
distribution by~\eqref{eqpropcond1}, whereas the second term
converges to $0$ in probability by~\eqref{eqneedL2}. This completes
the proof.
\end{pf}

\subsection{\texorpdfstring{Proofs of Theorems \protect\ref{theomainAss}--\protect\ref{theomainST}}
{Proofs of Theorems 1--3}}\label{sec2.4}
In this section we derive Theorems~\ref{theomainAss}--\ref{theomainST} as consequences of Proposition~\ref{propmain}.
We will replace condition~\eqref{eqpropcond2} by the following
one:
%
\begin{equation}\label{eqpropcond2a}
\sum_{\ome_1,\ome_2\in\Ome_n}r_{n}^2(\ome_1,\ome_2)=O(n s_n^2),\qquad
n\to\infty.
\end{equation}
Here, $r_n(t_1,t_2)=|t_1\cap t_2|$ is the overlap of the sets
$t_1,t_2\in T_n$. Condition~\eqref{eqpropcond2a} implies that~\eqref{eqpropcond2} holds with $\eps_n=1/\sqrt n$. Indeed, we have, as $n\to
\infty$,
\[
\sum_{\ome_1,\ome_2\in\Ome_n}\rho_{n}(\ome_1,\ome_2)1_{\rho_{n}(\ome
_1,\ome_2)>1/{\sqrt n} }
\leq
\sqrt n\sum_{\ome_1,\ome_2\in\Ome_n}\rho_{n}^2(\ome_1,\ome_2)
=o(s_n^2).
\]

\begin{pf*}{Proof of Theorem \protect\ref{theomainAss}}
The number of assignments on the set of $n$ elements is given by $|\Ass
_n|=n!$. To apply Proposition~\ref{propmain}, we take $E=\N\times\N$
and identify an assignment $\ome\in\Ome_n$ with the subset $\{(i,\ome
(i)); i=1,\ldots,n\}$ of~$E$.
To verify condition~\eqref{eqpropcond1} of Proposition~\ref{propmain}, note that
%
\begin{equation}\label{eqthasspr1}
\sum_{\ome\in\Ass_n}\XXX_n(\ome)
=\frac{1}{\sqrt n} \sum_{\ome\in\Ass_n} \sum_{i=1}^n \xi_{i,\ome(i)}
=\frac{(n-1)!}{\sqrt n}\sum_{i,j=1}^n \xi_{i,j}.
\end{equation}
It follows that $s_n^2$ defined in~\eqref{eqpropdef} is given by
%
\begin{equation}\label{eqthasspr2}
s_n^2=\Var\biggl[\sum_{\ome\in\Ass_n}\XXX_n(\ome)\biggr]=n!(n-1)!.
\end{equation}
The central limit theorem together with~\eqref{eqthasspr1} and~\eqref{eqthasspr2} implies that the random variable $s_n^{-1}\sum_{\ome\in
\Ass_n}\XXX_n(\ome)$ converges weakly to the standard Gaussian
distribution as $n\to\infty$. This verifies condition~\eqref{eqpropcond1} with $V\sim N(0,1)$.

Let us verify condition~\eqref{eqpropcond2a}. Let $\tilde\ome\in\Ome
_n$ be the identical assignment, that is, $\tilde\ome(i)=i$,
$i=1,\ldots,n$. We have
\[
\sum_{\ome_1,\ome_2\in\Ass_n} r_n^2(\ome_1,\ome_2)
=n!\sum_{\ome\in\Ass_n} r_n^2(\ome, \tilde\ome)
=n!\sum_{\ome\in\Ass_n} \Biggl(\sum_{i=1}^n 1_{\ome(i)=i}\Biggr)^2
=2(n!)^2,
\]
where the last equality follows from the well-known fact that the
expectation of the squared number of fixed points in a random
permutation is $2$. Together with~\eqref{eqthasspr2}, this verifies
condition~\eqref{eqpropcond2a}.
The proof is completed by applying Proposition~\ref{propmain}.
\end{pf*}

\begin{pf*}{Proof of Theorem~\protect\ref{theomainHP}}
To apply Proposition~\ref{propmain}, we take $E$ to be the set of all
two-element subsets of $\N$ and identify the set $V_n$ of vertices of
the complete graph $G_n$ with $\{1,\ldots,n\}$. Then, any (nonoriented)
Hamiltonian path $\ome\in\Ome_n$ can be viewed as a subset of $E$.
Let us verify condition~\eqref{eqpropcond1} of Proposition~\ref{propmain}.
The number of Hamiltonian paths in the complete graph~$G_n$, $n\geq3$, is given by $|\HP_n|=\frac12 (n-1)!$.
The number of Hamiltonian paths containing a given edge is easily seen
to be $(n-2)!$. Hence,
%
\begin{equation}\label{eqHPprws1}
\sum_{\ome\in\HP_n}\XXX_n(\ome)
=\frac{1}{\sqrt n}\sum_{\ome\in\HP_n} \sum_{e\in\ome} \xi_{e}
=\frac{1}{\sqrt n}\sum_{e\in E_n} \xi_{e} \sum_{\ome\in\HP_n} 1_{e\in
\ome}
=\frac{(n-2)!}{\sqrt n}\sum_{e\in E_n} \xi_{e}.
\end{equation}
Note that the number of edges in $G_n$ is $|E_n|=\frac12 n(n-1)$. It
follows that $s_n^2$ defined in~\eqref{eqpropdef} is given by
%
\begin{equation}\label{eqHPprws2}
s_n^2=\Var\biggl[\sum_{\ome\in\HP_n}\XXX_n(\ome)\biggr]=\frac12 (n-1)!(n-2)!.
\end{equation}
By the central limit theorem, combined with~\eqref{eqHPprws1}
and~\eqref{eqHPprws2}, the random variable $s_n^{-1}\sum_{\ome\in\HP
_n}\XXX_n(\ome)$ converges weakly to the standard Gaussian distribution
as $n\to\infty$. This verifies condition~\eqref{eqpropcond1} of
Proposition~\ref{propmain}.

We prove that~\eqref{eqpropcond2a} holds. We have
\begin{eqnarray}\label{eqHPprws3}
\sum_{\ome_1,\ome_2\in\HP_n}r_n^2(\ome_1,\ome_2)
&=&\sum_{\ome_1,\ome_2\in\HP_n}
\biggl(\sum_{e\in E_n} 1_{e\in\ome_1}1_{e\in\ome_2}\biggr)^2 \nonumber\\
&=&\sum_{\ome_1,\ome_2\in\HP_n}
\sum_{e,f\in E_n} 1_{e\in\ome_1}1_{e\in\ome_2}1_{f\in\ome_1}1_{f\in
\ome_2}\\
&=&
\sum_{e,f\in E_n}\biggl( \sum_{\ome\in\HP_n} 1_{e\in\ome}1_{f\in\ome
}\biggr)^2.\nonumber
\end{eqnarray}
The sum $\sum_{\ome\in\HP_n} 1_{e\in\ome}1_{f\in\ome}$ represents
the number of Hamiltonian paths containing the edges $e$ and $f$. If
$e=f$, then there are $(n-2)!$ such paths. If the edges $e$ and $f$
have exactly one common vertex, then the number of Hamiltonian paths
containing $e$ and $f$ is easily seen to be $(n-3)!$. Finally, if the
edges $e$ and $f$ do not have a common vertex, then the number of paths
containing both $e$ and $f$ is $2(n-3)!$. The number of pairs $(e,f)\in
E_n^2$ having exactly one common vertex is\vspace*{1pt} $6\bigl({{n}\atop{3}}\bigr)$, and the
number of pairs $(e,f)\in E_n^2$ without a common vertex is $6\bigl({{n}\atop{4}}\bigr)$. It follows from~\eqref{eqHPprws3} that $\sum_{\ome_1,\ome
_2\in\HP_n}r_n^2(\ome_1,\ome_2)$ is equal to
\[
\pmatrix{{n}\cr{2}}\bigl((n-2)!\bigr)^2+ 6\pmatrix{{n}\cr{3}}\bigl((n-3)!\bigr)^2+ 24\pmatrix{{n}\cr{4}}\bigl((n-3)!\bigr)^2.
\]
This expression is of order $O(ns_n^2)$ as $n\to\infty$. It follows
that~\eqref{eqpropcond2a} is fulfilled. The proof is completed by
applying Proposition~\ref{propmain}.
\end{pf*}
%
%
\begin{pf*}{Proof of Theorem~\protect\ref{theomainST}}
By Cayley's theorem, the number of spanning trees on $n$ vertices is
given by $|\ST_n|=n^{n-2}$. Since each spanning tree has $n-1$ edges,
and since there are $n(n-1)/2$ edges,\vadjust{\goodbreak} any edge is contained in
$2n^{n-3}$ trees. Hence,
%
\begin{equation}\label{eqSTwsp1}
\sum_{\ome\in\ST_n}\XXX_n(\ome)
=\frac1 {\sqrt{n-1}}\sum_{\ome\in\ST_n} \sum_{e\in\ome} \xi_{e}
=\frac{2n^{n-3}} {\sqrt{n-1}}\sum_{e\in E_n} \xi_{e}.
\end{equation}
It follows that
%
\begin{equation}\label{eqSTwsp2}
s_n^2=\Var\biggl[\sum_{\ome\in\ST_n}\XXX_n(\ome)\biggr]=2n^{2n-5}.
\end{equation}
By the central limit theorem together with~\eqref{eqSTwsp1} and~\eqref{eqSTwsp2}, the random variable $s_n^{-1}\sum_{\ome\in\ST_n}\XXX
_n(\ome)$ converges weakly to the standard Gaussian distribution as
$n\to\infty$.

Let us verify condition~\eqref{eqpropcond2a}. As in~\eqref{eqHPprws3}, we have
\[
\sum_{\ome_1,\ome_2\in\ST_n}r_n^2(\ome_1,\ome_2)
=
\sum_{e,f\in E_n}\biggl(\sum_{\ome\in\ST_n} 1_{e\in\ome}1_{f\in\ome
}\biggr)^2.
\]
Given two edges $e$ and $f$, we will compute the number of spanning
trees $N_n(e,f)=\sum_{\ome\in\Ome_n}1_{e\in\ome}1_{f\in\ome}$ in the
complete graph $G_n$ containing these two edges. For $e=f$, we have
shown that this number is equal to $2n^{n-3}$. We claim that if the
edges $e$ and $f$ have exactly one common vertex, then
$N_n(e,f)=3n^{n-4}$, whereas if $e$ and $f$ do not have common
vertices, then $N_n(e,f)=4n^{n-4}$. For completeness, we will prove
this by using the transfer current theorem giving an interpretation of
random spanning trees in terms of electric networks (see \cite{lyonsperesbook}, Section 8.2). It says that the probability that a
uniformly chosen spanning tree (in any finite graph) contains two given
edges $e$ and $f$ is given by the determinant
%
\begin{equation}\label{eqtransfercurrent}
\det
\pmatrix{Y(e,e) & Y(e,f)\vspace*{2pt}\cr
Y(f,e)& Y(f,f)}
,
\end{equation}
where $Y(g,h)$ denotes the (signed) current which flows through the
(somehow oriented) edge $h$ if a battery is hooked up between the ends
of the (somehow oriented) edge $g=(v_1,v_2)$ with such voltage that the
total current flowing through the graph is $1$. By Kirchhoff's laws and
symmetry reasons, we have $Y(g,g)=2/n$, $Y(g,h)=1/n$ if $h$ is of the
form $(v_1, v)$ for some vertex $v\neq v_2$, and $Y(g,h)=1/n$ if
$h=(v,v_2)$ for some vertex $v\neq v_1$. If $g$ and~$h$ have no
vertices in common, then $Y(g,h)=0$. Inserting this into~\eqref{eqtransfercurrent} and recalling that the total number of spanning
trees in $\Ome_n$ is $n^{n-2}$, we obtain the above mentioned formulae
for $N_n(e,f)$.

Recall from the proof of Theorem~\ref{theomainHP} that the number of
pairs $(e,f)\in E_n^2$ having exactly one common vertex\vspace*{1pt} is $6\bigl({{n}\atop{3}}\bigr)$, whereas the number of pairs $(e,f)\in E_n^2$ having no
vertices in common is $6\bigl({{n}\atop{4}}\bigr)$. Thus, 
\[
\sum_{\ome_1,\ome_2\in\Ass_n}r_n^2(\ome_1,\ome_2)=4\pmatrix{{n}\cr{2}}n^{2(n-3)}+
54\pmatrix{{n}\cr{3}} n^{2(n-4)}+ 96 \pmatrix{{n}\cr{4}}n^{2(n-4)}.
\]
The right-hand side is of order $O(ns_n^2)$ as $n\to\infty$. This
completes the proof of~\eqref{eqpropcond2a}.
\end{pf*}

\subsection{\texorpdfstring{Proof of Theorem \protect\ref{theomainPol}}{Proof of Theorem 4}}\label{sec2.5}
We will verify conditions~\eqref{eqpropcond1} and~\eqref{eqpropcond2} of Proposition~\ref{propmain}. Recall that $p_{k}(x)$
is the probability that a simple (nearest-neighbor) \mbox{$d$-dimensional}
random walk which starts at the origin, visits the site $x\in\Z^d$ at
time $k\in\N\cup\{0\}$.
Note that with $s_n$ defined by~\eqref{eqpropdef}, we have
\begin{eqnarray}
\sum_{\ome\in\Ome_n} \XXX_n(\ome)&=&\frac{(2d)^n}{\sqrt n}\sum
_{k=1}^n\sum_{x\in\Z^d} p_k(x)\xi_k(x),\label{eqthpolprwsp1a}\\
s_n^2&=&\frac{(2d)^{2n}}{n} \sum_{k=1}^{n}\sum_{x\in\Z^d}
p_{k}^2(x).\label{eqthpolprwsp1b}
\end{eqnarray}
First, we find an asymptotic formula for $\sum_{k=1}^{n}\sum_{x\in\Z^d}
p_{k}^2(x)$ as $n\to\infty$.
A~symmetry argument shows that $\sum_{x\in\Z^d} p_{k}^2(x) =
p_{2k}(0)$. Also, by the multidimensional local limit theorem (e.g.,
\cite{lawlerbook}, Section \ref{secsubgraphs}), $p_{2k}(0) \sim2^{1-d}(\pi
k/\break d)^{-d/2}$ as $k\to\infty$.
Thus, in the case $d\geq3$ we have
%
\begin{equation}\label{eqsumpsq}
S^2:=\sum_{k=1}^{\infty}\sum_{x\in\Z^d} p_{k}^2(x)<\infty.
\end{equation}
For $d=1,2$, we obtain the following asymptotics as $n\to\infty$:
%
\begin{equation}\label{eqwspom6}
\sum_{k=1}^{n}\sum_{x\in\Z^d} p_{k}^2(x)
\sim
2^{1-d}\sum_{k=1}^{n} \biggl(\frac{d}{\pi k}\biggr)^{d/2} 
\sim
\cases{
\displaystyle 2\sqrt{\frac{n}{\pi}}, & \quad$d=1,$\vspace*{2pt}\cr
\displaystyle\frac1 {\pi} \log n, &\quad$d=2.$
}
\end{equation}
%
In the case $d\geq3$, combining~\eqref{eqthpolprwsp1a}--\eqref{eqsumpsq}, we obtain the following
relation verifying condition~\eqref{eqpropcond1}:
\[
\frac1{s_n}\sum_{\ome\in\Pol_n}\XXX_n(\ome)
\todistr\frac1 S \sum_{k=1}^{\infty}\sum_{x\in\Z^d} p_{k}(x)\xi_{k}(x),
\]
where the series on the right-hand side converges in the $L^2$-sense.

In the case $d=1,2$, we will verify condition~\eqref{eqpropcond1} by
proving that the random variable $\frac{1}{s_n}\sum_{\ome\in\Ome_n}\XXX
_n(\ome)$ converges as $n\to\infty$ to the standard Gaussian distribution.
To this end, we will show that a triangular array in which the $n$th
row consists of the random variables $\{p_k(x)\xi_k(x); k=1,\ldots,n,
x\in\Z^d\}$ (with only finitely of them being nonzero) satisfies the
Lyapunov condition: for some $\delta>0$ and as $n\to\infty$,
%
\begin{equation}\label{eqpolprlyapunov}
\sum_{k=1}^{n} \sum_{x\in\Z^d} \E[|p_{k}(x)\xi_k(x)|^{2+\delta}]=o
\Biggl(\Biggl(\sum_{k=1}^{n}\sum_{x\in\Z^d} p_{k}^2(x)\Biggr)^{{(2+\delta)}/2}\Biggr).
\end{equation}
Note that $\sup_{x\in\Z^d} p_{k}(x)=O(k^{-d/2})$ as $k\to\infty$ by the
multidimensional local limit theorem (see \cite{lawlerbook}, Section
\ref{secsubgraphs}). Recalling the assumption $\E|\xi|^{2+\delta}<\infty$, we have
\begin{eqnarray*}
\sum_{k=1}^{n} \sum_{x\in\Z^d} \E[|p_{k}(x)\xi_k(x)|^{2+\delta}]
&=&
C\sum_{k=1}^{n} \sum_{x\in\Z^d} p_{k}^{2+\delta}(x)\\
&\leq&
C \sum_{k=1}^{n} \biggl(k^{-(1+\delta)d/2} \sum_{x\in\Z^d} p_{k}(x)
\biggr)\\
&=&
C \sum_{k=1}^{n} k^{-(1+\delta)d/2}\\
&\leq&
C n^{-(1+\delta)d/2+1}.
\end{eqnarray*}
It follows from~\eqref{eqwspom6} that for $d=1,2$, the Lyapunov
condition~\eqref{eqpolprlyapunov} holds. To complete the
verification of condition~\eqref{eqpropcond1} of Proposition~\ref{propmain},
recall~\eqref{eqthpolprwsp1a}, \eqref{eqthpolprwsp1b} and apply
the Lyapunov central limit theorem.

Let us verify condition~\eqref{eqpropcond2} for every $d\in\N$.
Arguing as in~\eqref{eqHPprws3}, we obtain
%
\begin{equation}
\sum_{\ome_1,\ome_2\in\Pol_n} r^2_{n}(\ome_1,\ome_2)
=
\mathop{\sum_{k_1,k_2=1,\ldots,n}}_{x_1,x_2\in\Z^d}
\biggl( \sum_{\ome\in\Pol_n} 1_{\ome(k_1)=x_1}1_{\ome(k_2)=x_2}\biggr)^2.
\end{equation}
The sum $\sum_{\ome\in\Pol_n} 1_{\ome(k_1)=x_1}1_{\ome(k_2)=x_2}$
counts the polymers $\ome\in\Ome_n$ with the property $\ome(k_1)=x_1$,
$\ome(k_2)=x_2$. For $1\leq k_1\leq k_2\leq n$, the number of such
paths is $(2d)^n p_{k_1}(x_1)p_{k_2-k_1}(x_2-x_1)$.
It follows that
\begin{eqnarray*}
\sum_{\ome_1,\ome_2\in\Pol_n} r^2_{n}(\ome_1,\ome_2)
&\leq&2\cdot(2d)^{2n}\sum_{1\leq k_1\leq k_2\leq n}\sum_{x_1,x_2\in\Z
^d}p_{k_1}^2(x_1)p_{k_2-k_1}^2(x_2-x_1)\\
&\leq&
2\cdot(2d)^{2n}\Biggl(\sum_{k=0}^n \sum_{x\in\Z^d} p_k^2(x)\Biggr)^2.
\end{eqnarray*}
%
With $\eps_n=n^{-1/4}$, it follows that for any dimension $d\in\N$,
\[
\mathop{\sum_{\ome_1,\ome_2\in\Pol_n}}_{r_{n}(\ome_1,\ome_2)> \eps_n
n}r_{n}(\ome_1,\ome_2)
\leq
\frac{1}{n^{3/4}}\sum_{\ome_1,\ome_2\in\Pol_n}r_{n}^2(\ome_1,\ome_2)
=o(n s_n^2),   \qquad  n\to\infty,
\]
%
where the last step follows from~\eqref{eqthpolprwsp1b} combined
with~\eqref{eqsumpsq} (in the case $d\geq3$) or~\eqref{eqwspom6}
(in the case $d=1,2$).
This verifies condition~\eqref{eqpropcond2}.
The proof of Theorem~\ref{theomainPol} can be now completed by
applying Proposition~\ref{propmain}. 

\subsection{\texorpdfstring{Proof of Theorem \protect\ref{theomainBRW}}{Proof of Theorem 5}}\label{sec2.6}
Given two vertices $t_1=(v_1,\ldots,v_n)\in\VVV$ and $t_2=(w_1,\ldots
,w_n)\in\VVV$ of length $n\in\N$ denote by\vadjust{\goodbreak} $r_n(\ome_1,\ome_2)=\min\{
i\in\N\dvtx  v_i\neq w_i\}-1$ the number of common ancestors, excluding
$\varnothing$, of $\ome_1$ and $\ome_2$. The next lemma will be needed in
the proof of Theorem~\ref{theomainBRW}.
\begin{lemma}\label{lemBRW}
Fix $k\in\N$. Define a stochastic process $\{V_n^{(k)}; n\in\N\}$ by
%
\begin{equation}\label{eqdefBRWsubmartingale}
V_n^{(k)}=\frac{1}{m^{2n}} \mathop{\sum_{\ome_1,\ome_2\in\Ome_n}}_{\ome
_1\neq\ome_2} r^k_n(\ome_1,\ome_2).
\end{equation}
Then, the limit $V^{(k)}_{\infty}:=\lim_{n\to\infty}V_n^{(k)}$ exists
in $(0,\infty)$ a.s. 
\end{lemma}
\begin{pf}
Let $\Aa_n=\sigma\{Z_t; l(t)<n\}$ be the $\sigma$-algebra generated by
the genealogical structure of the first $n$ generations of the
branching random walk. By definition, the random variable $V_n^{(k)}$
is $\Aa_n$-measurable. We will show that the sequence $\{V_n^{(k)}; n\in
\N\}$ is a submartingale with respect to the filtration $\{\Aa_n; n\in\N
\}$. 
We have
\[
V_{n+1}^{(k)}=\frac{1}{m^{2n+2}}\mathop{\sum_{\ome_1,\ome_2\in\Ome_n}}_{
\ome_1\neq\ome_2}Z_{\ome_1}Z_{\ome_2}r^k_n(\ome_1,\ome_2)+
\frac{1}{m^{2n+2}}\sum_{\ome\in\Ome_n}Z_{\ome}(Z_{\ome}-1)n^k.
\]
By our assumptions, $m=\E Z_t>1$ and $\gamma_2:=\E[Z_t(Z_t-1)]\in
(0,\infty)$. It follows that
%
\begin{equation}\label{eqprBRWsubmartingale}
\E\bigl[V_{n+1}^{(k)}|\Aa_n\bigr]=V_n^{(k)}+\frac{\gamma_2 n^k}{m^{2n+2}}|T_n|>
V_n^{(k)},
\end{equation}
whence the submartingale property. The sequence $\{V_n^{(k)}; k\in\N\}$
is bounded in $L^1$, since applying~\eqref{eqprBRWsubmartingale}
recursively, we obtain
%
\begin{equation}
\E\bigl[V_{n+1}^{(k)}\bigr]=\E\bigl[V_{n}^{(k)}\bigr]+\frac{\gamma_2n^k}{m^{n+2}}
=\cdots
=\gamma_2 \sum_{i=1}^{n}\frac{i^k}{m^{i+2}}.
\end{equation}
By the martingale convergence theorem, $V_{\infty}^{(k)}=\lim_{n\to
\infty}V_n^{(k)}$ exists in $[0,\infty)$ a.s. To see that the limit is
nonzero a.s., consider particles in generation $n$ which are offsprings
of some fixed particle in generation $1$. It is a classical fact that
the number of these offsprings divided by $m^{n-1}$ converges to an
a.s. nonzero random variable (see~\cite{harrisbook}, page 13). Since
for any of these two offsprings $\ome_1,\ome_2$, we have $r_n(\ome
_1,\ome_2)\geq1$, it follows that $V_{\infty}^{(k)}>0$ a.s.
\end{pf}
%
%
\begin{pf*}{Proof of Theorem~\protect\ref{theomainBRW}}
Given vertices $t_1,t_2\in\VVV$ of length $n\in\N$, note that
$
\rho_n(t_1,t_2):=\E[\XXX_n(t_1)\XXX_n(t_2)]=\frac1n r_n(t_1,t_2).
$
For $n\in\N$, let $s_n>0$ be a~random variable defined by
\[
s_n^2=\sum_{\ome_1,\ome_2\in\Ome_n}\rho_n(\ome_1,\ome_2).
\]
First we prove that we have a.s. finite random variables $V,W$ defined by
%
\begin{equation}\label{eqBRWcond1}
V=\lim_{n\to\infty}\frac{1}{s_n}\sum_{\ome\in\Ome_n}\XXX_n(\ome),\qquad
W=-\lim_{n\to\infty}\frac{\sqrt n}{|T_n|}\sum_{\ome\in\Ome_n}\XXX_n(\ome).\vspace*{-2pt}
\end{equation}
By Lemma~\ref{lemBRW}, we have
%
\begin{equation}\label{eqBRWasymptsn}
\qquad\lim_{n\to\infty}\sqrt n m^{-n} s_n =\lim_{n\to\infty} \sqrt
{V_n^{(1)}+m^{-2n}n|T_n|}=\sqrt{V_{\infty}^{(1)}}\in(0,\infty)\qquad
\mbox{a.s.},\vspace*{-2pt}
\end{equation}
where we have also used that $\lim_{n\to\infty}m^{-n}|T_n|$ exists in
$(0,\infty)$ a.s. (see \cite{harrisbook}, page 13).
It has been observed in~\cite{chen01} that $\{\sqrt nm^{-n}\sum_{\ome\in
\Ome_n}\XXX_n(\ome); n\in\N\}$ is an $L^2$-bounded martingale with
respect to the filtration $\{\Bb_n; n\in\N\}$, where $\Bb_n$ is the
$\sigma$-algebra generated by the genealogical structure $\{Z_t;
l(t)<n\}$ and the displacements $\{\xi_t; l(t)\leq n\}$ of the first
$n$ generations of the branching random walk. 
By the martingale convergence theorem and~\eqref{eqBRWasymptsn}, we
obtain that the limits in~\eqref{eqBRWcond1} exist a.s.
Also, it follows from Lemma~\ref{lemBRW} and~\eqref{eqBRWasymptsn} that
\begin{equation}
\lim_{n\to\infty}\frac1 {ns_n^2} \sum_{\ome_1,\ome_2\in\Ome
_n}r^2_n(\ome_1,\ome_2)
=
\lim_{n\to\infty}\frac{m^{2n}}{ns_n^2}\biggl( V_n^{(2)}+\frac
{n^2}{m^{2n}}|T_n|\biggr)
=\frac{V^{(2)}_{\infty}}{V_{\infty}^{(1)}},
\label{eqBRWcond2}\vspace*{-2pt}%
\end{equation}
which is finite a.s.

The proof of Theorem~\ref{theomainBRW} can be completed as follows.
Since the set $\Ome_n$ of particles in the $n$th generation is random,
we cannot apply Proposition~\ref{propmain} directly. To overcome this
difficulty, we will use a conditioning argument. We may assume that the
random variables $\{Z_t; t\in\VVV\}$ representing the numbers of
children are defined on a probability space $(\Omega_Z, \Aa_Z, \mu_Z)$
and the random variables $\{\xi_t; t\in\VVV\setminus\{\varnothing\}\}$
representing the displacements are defined on $(\Omega_{\xi},\Aa_{\xi},
\mu_{\xi})$. Then, we can define the branching random walk on the
product $\Omega_Z\times\Omega_{\xi}$ of both the spaces. Fix some $\tau
\in\Omega_Z$ and restrict all random variables to the set $\{\tau\}
\times\Omega_{\xi}$ endowed with the probability measure $\delta_{\tau
}\times\mu_{\xi}$, where $\delta_{\tau}$ is the Dirac measure at $\tau
$. Essentially, this means that we fix the realization of the
Galton--Watson tree but do not fix the displacements of the particles.
Note that the set $T_n$ becomes deterministic after such restriction.
It follows from~\eqref{eqBRWcond1} and~\eqref{eqBRWcond2} that
conditions~\eqref{eqpropcond1} and~\eqref{eqpropcond2a} of
Proposition~\ref{propmain} are fulfilled (in the restricted setting)
for $\mu_Z$-a.e. $\tau\in\Omega_Z$.
Applying Proposition~\ref{propmain}, we obtain that for $\mu_Z$-a.e. $\tau\in\Omega_Z$,
\[
\biggl\{\frac{\sqrt n}{|T_n|}\sum_{\ome\in\Ome_n} \bigl(1_{\XXX_n(\ome)\leq
z}- \Phi_n(z)\bigr); z\in\R\biggr\}
\stofd
\{p(z)W; z\in\R\},\vspace*{-2pt}
\]
where the random variables under consideration are restricted to the
space $\{\tau\}\times\Omega_{\xi}$. To complete the proof, integrate
over $\tau\in\Omega_Z$.
\end{pf*}

\section{\texorpdfstring{Proof of Theorem \protect\ref{theomainEA}}{Proof of Theorem 6}}

\vspace*{-3pt}\subsection{Hermite polynomials}\label{sechermitepoly}
$\!\!\!$We need to recall some facts about Hermite~po\-lynomials. Recall that
$p(z)=(2\pi)^{-1/2}e^{-z^2/2}$ is the standard\vadjust{\goodbreak} Gaussian densi\-ty. Let
$L^2(\R, p)$ be the set of all measurable functions\vspace*{-1pt} $f\dvtx \R\to\R$ such
that $\|f\|_{L^2(\R,p)}^2:=\int_{\R}f^2(z)p(z)\,dz$ is finite. The space
$L^2(\R, p)$ is a separable Hilbert space endowed with the scalar product
$\langle f,g \rangle_{L^2(\R, p)}=\int_{\R} f(z)g(z)p(z)\,dz$.
The (normalized)\vspace*{-1pt} Hermite polynomials $h_0, h_1,\ldots$ are defined\vspace*{-1pt} by
$h_n(z)=(-1)^n (n!)^{-1/2}e^{z^2/2}\frac{d^n}{dz^n}e^{-z^2/2}$. The
sequence $\{h_n\}_{n=0,1,\ldots}$ is an orthonormal basis in $L^2(\R
,p)$. For the proof of the next lemma see~\cite{sun65}, Lemma 1.1,
or~\cite{ivanovleonenkobook}, page 55.
\begin{lemma}\label{lemhermitegauss}
Let $(X,Y)$ be a zero-mean Gaussian vector with $\E X^2=\E Y^2=1$ and
$\E[XY]=\rho$. Then, for every $i,j\in\N\cup\{0\}$,
\[
\E[h_i(X)h_j(Y)]=
\cases{
\rho^i, &\quad$\mbox{if } i=j,$\vspace*{2pt}\cr
0, &\quad$\mbox{if } i\neq j.$
}
\]
\end{lemma}
%

Given $f\in L^2(\R, p)$ and $k\in\N$, we denote by $P_kf$ the
orthogonal projection of $f$ onto the orthogonal complement of the
$k$-dimensional linear subspace spanned by the first $k$ Hermite
polynomials $h_0,\ldots, h_{k-1}$. That is,
%
\begin{equation}\label{eqdefprojection}
(P_k f)(z)=\sum_{i=k}^{\infty} \langle f, h_i\rangle_{L^2(\R, p)}
h_i(z)=f(z)-\sum_{i=0}^{k-1} \langle f, h_i\rangle_{L^2(\R, p)} h_i(z).
\end{equation}
%

\begin{lemma}\label{lemhermite}
Let $(X,Y)$ be a zero-mean Gaussian vector with $\E X^2=\E Y^2=1$ and
$\E[XY]=\rho$. Then, for any $f,g\in L^2(\R, p)$ and $k\in\N$,
%
\begin{equation}
|\E[P_kf(X)P_kg(Y)]|\leq|\rho|^k \|f\|_{L^2(\R, p)}\|g\|_{L^2(\R, p)}.
\end{equation}
\end{lemma}
\begin{pf}
Write\vspace*{1pt} $f_i= \langle f, h_i\rangle_{L^2(\R,p)}$ and $g_i=\langle g,
h_i\rangle_{L^2(\R,p)}$ for $i\in\N\cup\{0\}$. We have $P_k f(X)=\sum
_{i=k}^{\infty} f_i h_i(X)$ and $P_k g(Y)=\sum_{i=k}^{\infty} g_i
h_i(Y)$. Using Lem\-ma~\ref{lemhermitegauss} and the inequality $|\rho
|\leq1$, we obtain
%
\[
|\E[P_kf(X)P_kg(Y)]|
=
\Biggl|\sum_{i=k}^{\infty} \rho^i f_ig_i\Biggr|
\leq|\rho|^k \sum_{i=0}^{\infty}|f_i||g_i|.
\]
To complete the proof, apply the Cauchy--Schwarz inequality.
\end{pf}
%

\subsection{Reduction method}\label{sec3.2}
The following proposition is an empirical central limit theorem for
Gaussian processes.
\begin{proposition}\label{propgauss2}
For every $n\in\N$, let $\{\XXX_n(\ome);\ome\in\Ome_n\}$ be a
zero-mean, unit-variance Gaussian process. Let $\rho_n(t_1,t_2)=\E[\XXX
_n(t_1)\XXX_n(t_2)]$ be the covariance function of $\XXX_n$. Define
$\varsigma_n\geq0$ by 
%
\begin{equation}
\varsigma_n^2:=\Var\biggl[\sum_{\ome\in\Ome_n}\bigl(\XXX_n^2(\ome)-1\bigr)
\biggr]=2\sum_{\ome_1,\ome_2\in\Ome_n} \rho_n^2(\ome_1,\ome_2).
\end{equation}
Suppose that for some random variable $V$ and for some sequence $\eps
_n>0$ satisfying $\lim_{n\to\infty}\eps_n=0$, the following three
conditions hold:
%
\begin{eqnarray}
\lim_{n\to\infty}\frac{1}{\varsigma_n^{2}}\sum_{\ome_1,\ome_2\in\Ome_n}
\rho_n(\ome_1,\ome_2)\hspace*{-5pt}&=&\hspace*{-5pt}0, \label{eqgausscond1}\\
\frac{1}{\varsigma_n}\sum_{\ome\in\Ome_n} \bigl(\XXX_n^2(\ome)-1\bigr)&\displaystyle\todistr&
V,\label{eqgausscond2}\\
\lim_{n\to\infty}\frac{1}{\varsigma_n^2}
\sum_{\ome_1,\ome_2\in\Ome_n} \rho_n^2(\ome_1,\ome_2)1_{|\rho_n(\ome
_1,\ome_2)|>\eps_n}\hspace*{-5pt}&=&\hspace*{-5pt}0.\label{eqgausscond3}
\end{eqnarray}
Then, the following convergence of stochastic processes holds true:
%
\begin{equation}
\biggl\{\frac1 {\varsigma_n} \sum_{\ome\in\Ome_n} \bigl(1_{\XXX_n(\ome)\leq
z}-\Phi(z)\bigr); z\in\R\biggr\} \tofd\biggl\{-\frac12 zp(z)V; z\in\R
\biggr\}.
\end{equation}
\end{proposition}
%
%
\begin{pf}
The proof is based on the reduction method of~\cite{taqqu74}. For
$x,z\in\R$, write $f(x;z)=1_{x\leq z}$. For $z\in\R$, define a
zero-mean random field $\{\YYY_n(\ome;z); \ome\in\Ome_n\}$ by
%
\begin{equation}
\YYY_n(\ome;z):=(P_3 f(\cdot; z))(\XXX_n(\ome)),
\end{equation}
where $P_3$ is the projection operator given in~\eqref{eqdefprojection}. Since the first three Hermite polynomials are
given by $h_0(x)=1$, $h_1(x)=x$, $h_2(x)=\frac1{\sqrt2} (x^2-1)$, this
means that
%
\begin{equation}\label{eqdefYnprop2}
\YYY_n(\ome;z)
=
1_{\XXX_n(\ome)\leq z}-\Phi(z)+p(z)\XXX_n(\ome)+\tfrac12 zp(z)\bigl(\XXX
_n^2(\ome)-1\bigr).
\end{equation}
By Lemma~\ref{lemhermite} with $k=3$ and $f=g$, we have
$
\E[\YYY_n(\ome_1;z)\YYY_n(\ome_2;z)]\leq C |\rho_n(\ome_1,\allowbreak\ome_2)|^3
$
for every $\ome_1,\ome_2\in\Ome_n$, where the constant $C$ does not
depend on $z\in\R$.
It follows that
\begin{eqnarray}\label{eqgaussprwsp1}
\Var\biggl[\sum_{\ome\in\Ome_n} \YYY_n(\ome;z)\biggr]
&=&\sum_{\ome_1,\ome_2\in\Ome_n} \E[\YYY_n(\ome_1;z)\YYY_n(\ome
_2;z)]\nonumber\\
&\leq& C \sum_{\ome_1,\ome_2\in\Ome_n} |\rho_n(\ome_1,\ome_2)|^3
\nonumber
\\[-8pt]
\\[-8pt]
\nonumber
&\leq& C\eps_n \sum_{\ome_1,\ome_2\in\Ome_n} \rho_n^2(\ome_1,\ome_2)+
C\mathop{\sum_{\ome_1,\ome_2\in\Ome_n}}_{|\rho_n(\ome_1,\ome_2)|> \eps
_n} \rho_n^2(\ome_1,\ome_2)\\
&=&o(\varsigma_n^2)\nonumber
\end{eqnarray}
as $n\to\infty$, where the last step follows from the assumption $\lim
_{n\to\infty}\eps_n=0$ and condition~\eqref{eqgausscond3}.
Take some $z_1,\ldots,z_d\in\R$. Then, for every $i=1,\ldots,d$, it
follows from~\eqref{eqdefYnprop2} that we have the following
decomposition:
\begin{eqnarray*}
&&\frac1 {\varsigma_n} \sum_{\ome\in\Ome_n} \bigl(1_{\XXX_n(\ome
)\leq z_i}-\Phi(z_i)\bigr)\\
&&\qquad=
-\frac{p(z_i)}{\varsigma_n}\sum_{\ome\in\Ome_n}\XXX_n(\ome)
-\frac{z_ip(z_i)}{2\varsigma_n}\sum_{\ome\in\Ome_n} \bigl(\XXX_n^2(\ome)-1\bigr)
\\
&&\qquad\quad{}+\frac{1}{\varsigma_n}\sum_{\ome\in\Ome_n} \YYY_n(\ome;z_i).
\end{eqnarray*}
As $n\to\infty$, the first term converges to $0$ in probability by
condition~\eqref{eqgausscond1}. The second term converges in
distribution to $-\frac12 z_ip(z_i)V$ by condition~\eqref{eqgausscond2}. Finally, the third term converges to $0$ in
probability by~\eqref{eqgaussprwsp1}. This completes the
proof.
\end{pf}

\subsection{\texorpdfstring{Completing the proof of Theorem \protect\ref{theomainEA}}
{Completing the proof of Theorem 6}}\label{sec3.3}
We will verify the conditions of Proposition~\ref{propgauss2}.
Given a spin configuration $\ome\in\Ome_n$ and an edge \mbox{$e=\{v_1,v_2\}
\in E_n$}, we write $\ome\diamond e=\ome(v_1)\ome(v_2)\in\{+1,-1\}$.
Recall that the energy of a spin configuration $\ome\in\Ome_n$ is given by
%
\begin{equation}\label{eqdefEA1}
\XXX_n(\ome)=|E_n|^{-1/2} \sum_{e\in E_n}(\ome\diamond e) J(e),
\end{equation}
where $\{J(e); e\in E_n\}$ are independent standard Gaussian random variables.

We start by verifying condition~\eqref{eqgausscond1} of
Proposition~\ref{propgauss2}. Since\break \mbox{$\sum_{\ome\in\Ome_n}(\ome
\diamond e)=0$} for every edge $e\in E_n$, we have
\begin{eqnarray*}
\sum_{\ome\in\Ome_n} \XXX_n(\ome)
&=&
|E_n|^{-1/2}\sum_{\ome\in\Ome_n} \sum_{e\in E_n} (\ome\diamond e) J(e)\\
&=&
|E_n|^{-1/2}\sum_{e\in E_n}J(e) \sum_{\ome\in\Ome_n}(\ome\diamond e)\\
&=&0.
\end{eqnarray*}
Hence, $\sum_{\ome_1,\ome_2\in\Ome_n}\rho_n(\ome_1,\ome_2)=0$, which
implies that condition~\eqref{eqgausscond1} holds.

Let us verify condition~\eqref{eqgausscond2} of Proposition~\ref{propgauss2}. Note that for every different edges $e_1,e_2\in E_n$,
we have $\sum_{\ome\in\Ome_n}(\ome\diamond e_1)(\ome\diamond e_2)=0$. Hence,
\begin{eqnarray}\label{eqgausswspom1}
&&\sum_{\ome\in\Ome_n} \bigl(\XXX_n^2(\ome)-1\bigr)\nonumber\\
&&\qquad=
|E_n|^{-1}\sum_{\ome\in\Ome_n} \sum_{e_1,e_2\in E_n} \bigl((\ome\diamond
e_1)(\ome\diamond e_2) J(e_1)J(e_2)-1_{e_1=e_2}\bigr)\\
&&\qquad=
\frac{|\Ome_n|}{|E_n|} \sum_{e\in E_n} \bigl(J^2(e)-1\bigr).\nonumber
\end{eqnarray}
By Lemma~\ref{lemhermitegauss}, $\E[(\XXX_n^2(\ome_1)-1)(\XXX_n^2(\ome
_2)-1)]=2\rho^2_n(\ome_1,\ome_2)$. It follows from this and~\eqref{eqgausswspom1} that
%
\begin{equation}\label{eqgausswspom2}
\varsigma_n^2
=
\Var\biggl[\sum_{\ome\in\Ome_n}\bigl(\XXX_n^2(\ome)-1\bigr)\biggr]
=
\frac{|\Ome_n|^2}{|E_n|^2}\Var\biggl[ \sum_{e\in E_n} \bigl(J^2(e)-1\bigr)\biggr]
=
\frac{2|\Ome_n|^2}{|E_n|}.
\end{equation}
The central limit theorem together with~\eqref{eqgausswspom1}
and~\eqref{eqgausswspom2} implies that condition~\eqref{eqgausscond2} is satisfied with $V\sim N(0,1)$.

Let us verify condition~\eqref{eqgausscond3} of Proposition~\ref{propgauss2}. It follows from~\eqref{eqdefEA1} that for every $\ome
_1,\ome_2\in\Ome_n$, $\rho_n(\ome_1,\ome_2)=|E_n|^{-1}\sum_{e\in
E_n}(\ome_1\diamond e)(\ome_2\diamond e)$. Define a spin configuration
$\tilde\ome\in\Ome_n$ by requiring that $\tilde\ome(v)=1$ for every
vertex $v\in V_n$. It follows that
\begin{eqnarray}\label{eqprsp1}
\sum_{\ome_1,\ome_2 \in\Ome_n} \rho_n^4(\ome_1,\ome_2)
&=&
|\Ome_n|\sum_{\ome\in\Ome_n} \rho_n^4(\tilde\ome,\ome)\nonumber\\
&=&
|\Ome_n||E_n|^{-4}\sum_{\ome\in\Ome_n} \biggl(\sum_{e\in E_n}(\ome
\diamond e)\biggr)^4\\
&=&
|\Ome_n||E_n|^{-4} \sum_{e_1,e_2,e_3,e_4\in E_n} \sum_{\ome\in\Ome_n}
\prod_{k=1}^4(\ome\diamond e_k).\nonumber
\end{eqnarray}
It will be convenient to write $\eta(e_1,\ldots,e_4)=\sum_{\ome\in\Ome
_n}\prod_{k=1}^4 (\ome\diamond e_k)$. If some vertex $v\in V_n$ belongs
to exactly one or exactly three of the edges $e_1,\ldots,e_4$, then
$\eta(e_1,\ldots,e_4)=0$ by spin flip symmetry. Consider some quadruple
$e_1,\ldots,e_4$ for which $\eta(e_1,\ldots,e_4)\neq0$. We will show
that there are at most $C|E_n|^2$ such quadruples. The union of all
vertices belonging to $e_1,\ldots,e_4$ consists of $2$ or $4$ elements.
In both cases, we can find $i,j\in1,\ldots,4$ such that the union of
vertices belonging to $e_1,\ldots,e_4$ coincides with the union of the
vertices of~$e_i,e_j$. There are at most $|E_n|^2$ possibilities to
choose $e_i$ and $e_j$ and a bounded number of choices for the
remaining two edges. To summarize, there are at most $C|E_n|^2$ terms
of the form $\eta(e_1,\ldots,e_4)$ which are nonzero, and any such term
is bounded by $|\Ome_n|$. It follows from these considerations
and~\eqref{eqprsp1} that
%
\begin{equation}\label{eqgausswspom3}
\sum_{\ome_1,\ome_2 \in\Ome_n} \rho_n^4(\ome_1,\ome_2)
\leq|\Ome_n||E_n|^{-4}\cdot C |E_n|^2|\Ome_n|
\leq C \frac{|\Ome_n|^2}{|E_n|^{2}}.
\end{equation}
Now we are able to verify condition~\eqref{eqgausscond3}. Since $\lim
_{n\to\infty}|E_n|=\infty$, we can choose $\eps_n>0$ in such a way that
$\lim_{n\to\infty}\eps_n=0$ but $\lim_{n\to\infty} \eps_n^2|E_n|=\infty
$. Recalling~\eqref{eqgausswspom2} and~\eqref{eqgausswspom3}, we obtain
\[
\frac{1}{\varsigma_n^2}
\sum_{\ome_1,\ome_2\in\Ome_n} \rho_n^2(\ome_1,\ome_2)1_{|\rho_n(\ome
_1,\ome_2)|>\eps_n}
\leq
\frac{1}{\varsigma_n^2 \eps_n^2} \sum_{\ome_1,\ome_2\in\Ome_n} \rho
_n^4(\ome_1,\ome_2)
\leq
\frac{C}{\eps_n^2|E_n|},
\]
which converges to $0$ as $n\to\infty$. This completes the verification
of condition~\eqref{eqgausscond3} of Proposition~\ref{propgauss2}.

\section*{Acknowledgments}
The author expresses his gratitude to Wolfgang Kar\-cher, Daniel
Meschenmoser and Florian Timmermann for useful discussions on empirical
central limit theorems.


\printaddresses

\end{document}